\documentclass{article}

\usepackage{amsfonts, latexsym, amssymb,amsmath,amstext}
\usepackage[dvips]{graphics,color}

\title{Transformations of
hypergeometric elliptic integrals}

\author{Raimundas Vid\=unas\\
\em Kobe University}

\newtheorem{theorem}{Theorem}[section]

\usepackage{latexsym}

\newcommand{\exppi}[1]{{\bf e}^{2\pi i #1}}
\newcommand{\hpg}[5]{{}_{#1}\mbox{\rm F}_{\!#2}\!
  \left(\left.{#3 \atop #4}\right| #5 \right) }

\newcommand{\proof}{{\bf Proof. }}
\newcommand{\qed}{\hfill $\Box$}
\newcommand{\equal}{\!\!\!=\!\!\!}
\newcommand{\CC}{{\mathbb C}}

\newcommand{\ZZ}{\mbox{\bf Z}}

\date{}

\begin{document}

\maketitle

\begin{abstract} 
The paper classifies algebraic transformations of
Gauss hypergeometric functions with the local exponent differences
$(1/2,1/4,1/4)$, $(1/2,1/3,1/6)$ and $(1/3,1/3,1/3)$. 
These form a special class of algebraic transformations of Gauss hypergeometric functions,
of arbitrary high degree. The Gauss hypergeometric functions can be identified as
elliptic integrals on the genus 1 curves $y=x^3-x$ or $y=x^3-1$. 
Especially interesting are algebraic transformations of the hypergeometric functions into themselves;
these transformations come from isogenies of the respective elliptic curves.
\end{abstract}

\section{Introduction}
\label{ellints}

Well-known quadratic transformations of Gauss hypergeometric functions
were supplemented by Goursat \cite{goursat} with transformations of degree 3, 4 and 6.
It was widely assumed that there are no other algebraic transformations, 
except between algebraic hypergeometric functions; see \cite[Section 2.1.5]{bateman} for example.
However, this assumption is wrong, as noticed in  \cite{algtgauss} and \cite{andkitaev}.
All possible algebraic transformations coming from pull-back transformations of their
second order Fuchsian differential equations are reviewed in \cite{algtgauss}.
Particularly interesting cases are transformations of Gauss hypergeometric functions
with the local exponent differences $(1/k,1/\ell,1/m)$ such that $k,\ell,m$ are positive integers.
If $1/k+1/\ell+1/m>1$, the transformed hypergeometric functions are algebraic.
By Klein's theorem, any algebraic Gauss hypergeometric function is a pull-back transformation
of a standard Gauss hypergeometric function with the local exponent differences 
satisfying $1/k+1/\ell+1/m>1$;
computation of these transformations is thoroughly considered in \cite{talggaus}. 
There is a finite list of transformations for Gauss hypergeometric functions with the local exponent differences satisfying $1/k+1/\ell+1/m<1$, of degree up to 24; see \cite{thyperbolic}.

This paper considers transformations of Gauss hypergeometric functions 
with the local exponent differences $(1/k,1/\ell,1/m)$ satisfying $1/k+1/\ell+1/m=1$, 
where $k,\ell,m$ are positive integers. Up to permutations, there are three such
triples of local exponent differences: $(1/2,1/4,1/4)$,
$(1/2,1/3,1/6)$, and $(1/3,1/3,1/3)$. Euler's hypergeometric equations
with these local exponent differences have trivial solutions.
The degenerate structure of the 24 Kummer's solutions is explained
in Appendix Section \ref{appendix}.

We particularly consider the following hypergeometric functions and their
integral representations:
\begin{eqnarray} \label{elin4}
\hpg{2}{1}{1/2,\,1/4}{5/4}{z} & = & \frac{z^{-1/4}}{4}\, \int_0^z
t^{-3/4}\,(1-t)^{-1/2}\,dt,\\ \label{elin6}
\hpg{2}{1}{1/2,\,1/6}{7/6}{z} & = & \frac{z^{-1/6}}{6}\, \int_0^z
t^{-5/6}\,(1-t)^{-1/2}\,dt,\\ \label{elin3}
\hpg{2}{1}{1/3,\,2/3}{4/3}{z} & = & \frac{z^{-1/3}}{3}\, \int_0^z
t^{-2/3}\,(1-t)^{-2/3}\,dt.
\end{eqnarray}
The local exponent differences of these Gauss hypergeometric functions are,
respectively, $(1/2,1/4,1/4)$, $(1/2,1/3,1/6)$, and $(1/3,1/3,1/3)$. 
The paths of integration should lie in a sector that contains $z$,
and the branch of fractional power of $t$ should correspond to the
branch of the fractional power of $z$ outside the integral.
Expressions and transformations for other hypergeometric functions 
with the same local exponent differences follow from
connection formulas in Appendix Section \ref{appendix}. 

As we shall see, the three integrals above are elliptic integrals,
as algebraic relations between the integrands and the
variable $t$ define curves of genus 1. The first integral is
defined on the curve isomorphic to $y^2=x^3-x$. The other two integrals 
are easily transformed to elliptic integrals on the isomorphic curves
$y^2=x^3-1$ and $X^3+Y^3=1$, respectively. 

A general algebraic transformation of Gauss hypergeometric functions is an identity
of the form
\begin{equation} \label{hpgtransf}
\hpg{2}{1}{\!\widetilde{A},\,\widetilde{B}\,}{\widetilde{C}}{\,x}
=\theta(x)\;\hpg{2}{1}{\!A,\,B\,}{C}{\varphi(x)},
\end{equation}
where $\varphi(x)$ is a rational function of $x$, and $\theta(x)$ is a
radical function, i.e., product of some powers of rational functions.
By the general classification scheme of algebraic transformations in \cite{algtgauss},
algebraic transformations of the functions in (\ref{elin4})--(\ref{elin3})
should satisfy the following properties:
\begin{itemize}
\item The covering can ramify only above the three singularities of the differential equation
under the corresponding pull-back transformation.
\item There should be precisely $d+2$ distinct points above the
three singularities of the original equation (by Hurwitz formula).
\item There should be precisely $3$ singular points (of the transformed equation)
above the three singularities of the original equation.
\item The positive local exponent differences of the transformed equation
should sum up to $1$ as well. This follows from the degree formula \cite[(13)]{algtgauss}.
\end{itemize}
The possible ramification patterns are shown in Table \ref{elltab}.
\begin{table}
\begin{center} \begin{tabular}{|c|c|c|c|}
\hline \multicolumn{2}{|c|}{Local exponent differences} & Degree &
Ramification above \\ \cline{1-2}
$\!(1/k,\,1/\ell,\,1/m)\!$ & above & $d$ & the regular singular points \\
 \hline
$(1/2,\,1/4,\,1/4)$ & $(1/2,\,1/4,\,1/4)$ & $4n$ & $2n\!*2=n\!*4=(n\!-\!1)\!*4+2+1+1$ \\
$(1/2,\,1/4,\,1/4)$ & $(1/2,\,1/4,\,1/4)$ & $4n\!+\!1$ & $2n\!*2+1=n\!*4+1=n\!*4+1$ \\
$(1/2,\,1/4,\,1/4)$ & $(1/2,\,1/4,\,1/4)$ & $4n\!+\!2$ & $(2n\!+\!1)\!*2=n\!*4+2=n\!*4+1+1$ \\
$(1/2,\,1/3,\,1/6)$ & $(1/2,\,1/3,\,1/6)$ & $6n$ & $3n\!*2=2n\!*3=(n\!-\!1)\!*6+3+2+1$ \\
$(1/2,\,1/3,\,1/6)$ & $(1/2,\,1/3,\,1/6)$ & $6n\!+\!1$ & $3n\!*2+1=2n\!*3+1=n\!*6+1$ \\
$(1/2,\,1/3,\,1/6)$ & $(1/2,\,1/3,\,1/6)$ & $6n\!+\!3$ &
$(3n\!+\!1)\!*2\!+\!1=(2n\!+\!1)\!*3=n\!*6\!+\!2\!+\!1$ \\
$(1/2,\,1/3,\,1/6)$ & $(1/2,\,1/3,\,1/6)$ & $6n\!+\!4$ &
$(3n\!+\!2)\!*2=(2n\!+\!1)\!*3\!+\!1=n\!*6\!+\!3\!+\!1$ \\
$(1/2,\,1/3,\,1/6)$ & $(1/3,\,1/3,\,1/3)$ & $6n$ & $3n\!*2=2n\!*3=(n\!-\!1)\!*6+2+2+2$ \\
$(1/2,\,1/3,\,1/6)$ & $(1/3,\,1/3,\,1/3)$ & $6n$ & $3n\!*2=(2n\!-\!1)\!*3+1+1+1=n\!*6$ \\
$(1/2,\,1/3,\,1/6)$ & $(1/3,\,1/3,\,1/3)$ & $6n\!+\!2$ & $(3n\!+\!1)\!*2=2n\!*3+1+1=n\!*6+2$ \\
$(1/2,\,1/3,\,1/6)$ & $(1/3,\,1/3,\,1/3)$ & $6n\!+\!4$ &
$(3n\!+\!2)\!*2=(2n\!+\!1)\!*3\!+\!1=n\!*6\!+\!2\!+\!2$ \\
$(1/2,\,1/3,\,1/6)$ & $(2/3,\,1/6,\,1/6)$ & $6n$ & $3n\!*2=2n\!*3=(n\!-\!1)\!*6+4+1+1$ \\
$(1/2,\,1/3,\,1/6)$ & $(2/3,\,1/6,\,1/6)$ & $6n\!+\!2$ & $(3n\!+\!1)\!*2=2n\!*3+2=n\!*6+1+1$ \\
$(1/3,\,1/3,\,1/3)$ & $(1/3,\,1/3,\,1/3)$ &
$3n$ & $n\!*3=n\!*3=(n\!-\!1)\!*3+1+1+1$ \\
$(1/3,\,1/3,\,1/3)$ & $(1/3,\,1/3,\,1/3)$ & $3n\!+\!1$ &
$n\!*3\!+\!1=n\!*3\!+\!1=n\!*3\!+\!1$ \\ \hline
\end{tabular} \end{center}
\caption{Transformations of hypergeometric elliptic integrals}
\label{elltab}
\end{table}
In the last column, multiplicative terms give a ramification order after the * sign, and 
the number of points with that ramification order in front; 
the = signs separate ramification orders in the three fibers.

Most of the entries in Table \ref{elltab} represent transformations of hypergeometric equations
with the local exponent differences $(1/k,1/\ell,1/m)$ to themselves. Consequently, 
the hypergeometric functions in (\ref{elin4})--(\ref{elin3}) can be transformed to themselves,
giving formulas (\ref{hpgtransf}) with the same arguments $A,B,C$ on both sides.
The main result of this paper is that these algebraic transformations are induced by
isogeny endomorphisms of the respective elliptic curves. This proposition was
immediately suggested by Frits Beukers in a private correspondence \cite{beukrspriv}.
The integrands in (\ref{elin4})--(\ref{elin3}) are holomorphic differentials 
on the respective elliptic curves; an isogeny transforms it to a holomorphic differential
again, necessarily proportional to the original one. This gives a transformation
of the hypergeometric function into itself.

Transformations from the local exponent differences $(1/2,\,1/3,\,1/6)$ to
$(1/3,\,1/3,\,1/3)$ or $(2/3,\,1/6,\,1/6)$ turn out to be compositions of 
fixed quadratic transformations and the mentioned isogeny transformations.
This is considered in Section \ref{otherells}. A particular implication is that
there are no pull-back transformations from $(1/2,1/3,1/6)$  to $(1/3,1/3,1/3)$
of degree $6n+4$, as there are no isogeny transformations of degree $3n+2$
for either of the two hypergeometric equation. A Gauss hypergeometic function
with the local exponent differences $(2/3,\,1/6,\,1/6)$ is a hyperelliptic integral
on a genus 2 curve isomorphic to $Y^2=X^6+1$. 

Except in Section \ref{otherells}, by an isogeny we mean an isogeny endomorphism
on an elliptic curve. Recall that an isogeny is an algebraic map between elliptic curves
that respects the addition law on them. Isogeny endomorphisms on an elliptic curve
form a ring; addition of isogenies is given by the addition law on the elliptic curve,
while multiplication of isogenies is their composition.


\section{Elliptic integrals on the curve $y^2=x^3-x$}

Let $E_1$ denote the genus 1 curve defined by the equation
$y^2=x^3-x$. We consider $E_1$
as an elliptic curve in the usual way --- by specifying the point
at infinity as the origin of $E_1$ as an additive group. 

Let $f_4=t^{-3/4}\,(1-t)^{-1/2}$ denote the integrand function
in (\ref{elin4}). The algebraic curve with the function field
$\CC(t,f_4)$ is isomorphic to $E_1$ by the isomorphisms
$(t,f_4)=\left(x^{-2},\,x^3y^{-1}\right)$ and $x=f_4^2\,t\,(1-t)$.
Therefore (\ref{elin4}) is an elliptic integral on $E_1$. 
By substituting $t=x^{-2}$ in the integral 
(\ref{elin4}) we get:
\begin{equation} \label{elin4a}
\hpg{2}{1}{1/2,\,1/4\,}{5/4}{\,z}=\frac{z^{-1/4}}{2}\;
\int_{1/\!\sqrt{z}}^{\infty}\;\frac{dx}{\sqrt{x^3-x}}.
\end{equation}
The most convenient path of integration is the constant phase
path along the ray from $1/\sqrt{z}$ to the infinity. 

In the next two theorems we show that transformations of the hypergeometric function 
in (\ref{elin4a}) into itself correspond to isogeny endomorphisms on $E_1$.
As $E_1$ has {\em complex multiplication} by $i$,
the ring of isogenies on $E_1$ is isomorphic to the ring $\ZZ[i]$
of Gaussian integers; see \cite{silverman1}. We identify
$i\in\ZZ[i]$ with the automorphism $(x,y)\mapsto (-x,iy)$. For
$a+bi\in\ZZ[i]$, let $\phi_{a+bi}=(\varphi_{a+bi},\psi_{a+bi})$
denote the corresponding isogeny; its degree is equal
to the norm $a^2+b^2$ of $a+bi$. Here are a few isogenies on $E_1$:
\begin{eqnarray} \label{isog1}
\phi_{1+i}: & (x,y)\,\mapsto\!\! & \left(\frac{x^2-1}{2i\,x},
\;\frac{y\;(x^2+1)}{2(i-1)\,x^2}\right),\\
\phi_{2}: & (x,y)\,\mapsto\!\! &\left(\frac{(x^2+1)^2}{4\,x\,(x^2-1)},
\frac{(x^2+1)(x^4-6x^2+1)}{8\,x\,y\,(x^2-1)}\right),\\
\label{isog3} \phi_{1+2i}: & (x,y)\,\mapsto\!\! & \left(
\frac{x(x^2\!-\!1\!-\!2i)^2}{((1\!+\!2i)x^2\!-\!1)^2},
\frac{y(x^4\!+\!(2\!+\!8i)x^2\!+\!1)(x^2\!-\!1\!-\!2i)}
{((1\!+\!2i)x^2\!-\!1)^3}\right).
\end{eqnarray}

\begin{theorem} \label{th:el4tr}
Suppose that $\phi:E_1\to E_1$ is a non-zero isogeny on $E_1$, of degree $d$.
The isogeny transforms the $x$-coordinate as $x \mapsto x\,\mu(x^2)$,
where $\mu(s)$ is a rational function with a finite non-zero limit $\mu_0=\lim_{s\to\infty}\mu(s)$. 
We have the identity
\begin{equation} \label{el4tr}
\hpg{2}{1}{1/2,\,1/4}{5/4}{\,z} = \left(\frac{\mu(1/z)}{\mu_0}\right)^{-1/2}\;
\hpg{2}{1}{1/2,\,1/4}{5/4}{\frac{z}{\mu\!\left(1/z\right)^2}}.
\end{equation}
The degree of the rational function $z/\mu\!\left(1/z\right)^2$ is equal to $d$.
\end{theorem}
\proof Let $\varphi(x)=\varphi(x,y)$ and $\psi(x)=\psi(x,y)$ denote, respectively, 
the $x,y$ components of the isogeny $\phi$, so that 
$\phi(x,y)=\left(\varphi(x,y),\psi(x,y)\right)$.

The integrand in (\ref{elin4a}) is a holomorphic differential $1$-form on $E_1$.
The substitution $x\mapsto \varphi(x,y)$ into the integral in (\ref{elin4a}) must be
an integral of  a holomorphic differential $1$-form on $E_1$ again. 
Since the linear space of holomorphic differentials on an elliptic curve is one-dimensional,
the transformed differential form must be proportional to $dx/\sqrt{x^3-x}$. 
Moreover, the upper integration bound does not change, because isogenies fix the point at infinity. 
The lower integration bound is transformed as $z\mapsto\varphi(1/\sqrt{z})^{-2}$.
Once we show that $\varphi(x)=x\,\mu(x^2)$ as stated in theorem, transformation (\ref{el4tr})
follows; the power factor must evaluate to 1 at $z=0$.

The addition law on $E_1$ gives addition of isogenies. For $u,v\in\ZZ[i]$ we have:
\begin{equation} \label{elladd}
\varphi_{u+v}=\lambda^2-\varphi_u-\varphi_v, \qquad
\psi_{u+v}=-\psi_u-\lambda\,(\varphi_{u+v}-\varphi_u)
\end{equation}
where
\begin{equation}
\lambda=\frac{\psi_u-\psi_v}{\varphi_u-\varphi_v}\qquad\mbox{or}\qquad
\lambda=\frac{\varphi_u^2+\varphi_u\varphi_v+\varphi_v^2-1}{\psi_u+\psi_v}.
\end{equation}
By induction on the norm $|u|$ we conclude that any isogeny $\phi_u$ can be written
in the form $\left(x\,\mu(x^2),y\,\nu(x^2)\right)$, where $\mu(s)$, $\nu(s)$ are 
rational functions satisfying $\lim_{s\to\infty} \mu(s)=1/u^2$ and
$\lim_{s\to\infty} \nu(s)=1/u^3$. Transformation (\ref{el4tr}) follows.

For the last statement, note that the point $(x,y)=(0,0)$ is in the kernel
of all isogenies $\phi_u$ with even $d=|u|$; hence $s=0$ is a pole of $\mu(s)$ when $d$ is even.
The degree of the rational function $\varphi(x)=x\,\mu(x^2)$ is $d$, and the degree
of $\mu(s)$ is equal to $\lfloor d/2\rfloor$ for both even and odd $d$.
The factor $z$ in $z/\mu(1/z)^2$ simplifies in the same way, exactly when $d$ is even.
\qed\\

\noindent
Here are the transformations induced by isogenies (\ref{isog1})--(\ref{isog3}):
\begin{eqnarray}
\hpg{2}{1}{1/2,\,1/4}{5/4}{z} & = & \frac{1}{\sqrt{1-z}}\;
\hpg{2}{1}{1/2,\,1/4}{5/4}{-\frac{4\,z}{(z-1)^2}},\\
\hpg{2}{1}{1/2,\,1/4}{5/4}{z} & = & \frac{\sqrt{1-z}}{1+z}\;
\hpg{2}{1}{1/2,\,1/4}{5/4}{\frac{16\,z\,(z-1)^2}{(z+1)^4}},\\
\hpg{2}{1}{1/2,\,1/4}{5/4}{z}&=&\frac{1-z/(1\!+\!2i)}{1-(1\!+\!2i)z}\;
\hpg{2}{1}{1/2,\,1/4}{5/4}{\frac{z\,(z-1-2i)^4}{\big((1\!+\!2i)z-1\big)^4}}.
\end{eqnarray}
The first two formulas are special cases of classical algebraic transformations.
Transformations of the hypergeometric equation for (\ref{elin4a})
into itself form a group under the composition. Since the
isogenies $\varphi_{-1}$, $\varphi_i$, and $\varphi_{-i}$ induce
the trivial transformation of (\ref{elin4a}),  the group of the
hypergeometric transformations is isomorphic to the multiplicative
group $\ZZ[i]^*/\{\pm 1,\pm i\}$.

The following theorem implies that a transformation (\ref{el4tr}) of degree $d$ 
exists only if there is an element of $a+bi\in\ZZ[i]$ with the norm $d=a^2+b^2$.
In particular,  there is no transformation of degree 21, even if $21\mbox{ mod }4=1$,
as was hinted in \cite{beukrspriv}. 

\begin{theorem}
Any transformation of the form $(\ref{el4tr})$ comes from an isogeny endomorphism of $E_1$
as described in Theorem $\ref{th:el4tr}$.
\end{theorem}
\proof The first three entries in Table \ref{elltab} give explicit ramification patterns of 
possible pull-back coverings of Euler's hypergeometric equation with the local exponent
differences $(1/2,1/4,1/4)$ into itself. For identity $(\ref{el4tr})$ we should have the local exponent
difference $1/4$ at the points $z=0$ and $z=\infty$ on both projective lines above and below,
and the points $z=0$ must lie above each other. 


With this setting, a pull-back covering of degree $4n+1$ has the form:
\begin{equation}
z\mapsto \frac{z\,P(z)^4}{Q(z)^4}, \qquad 1-z\mapsto
\frac{(1-z)\,R(z)^2}{Q(z)^4},
\end{equation}
where $P(z)$, $Q(z)$ and $R(z)$ are polynomials of degree
$n$, $n$ and $2n$, respectively. The polynomial $Q(z)$ can be assumed to be monic.
The above form of the transformation implies the polynomial identity
\begin{equation} \label{pqrcheck1}
(1-z)\,R(z)^2=Q(z)^4-z\,P(z)^4.
\end{equation}
Similarly, pull-back coverings of degree $4n+2$ and $4n$ give the following polynomial identities,
respectively,
\begin{eqnarray} \label{pqr4check2}
R(z)^2 &\equal& (1-z)^2\,Q(z)^4-z\,P(z)^4,\\ \label{pqr4check3}
R(z)^2 &\equal& Q(z)^4-z\,(1-z)^2\,P(z)^4.
\end{eqnarray}
In the first identity, the polynomials $P(z)$, $Q(z)$ and $R(z)$ have degree 
$n$, $n$ and $2n+1$, respectively, while in the second identity the degrees 
are $n-1$, $n$ and $2n$.
 
In the three cases, the following maps can be checked to be isogeny endomorphisms
of the elliptic curve $E_1$, respectively:
\begin{eqnarray} \label{el4isog1}
(x,y)&\!\mapsto\!&\left(\frac{x\,Q(x^{-2})^2}{P(x^{-2})^2},\;
\frac{y\;Q(x^{-2})\,R(x^{-2})}{P(x^{-2})^3}\right),\\ \label{el4isog2}
(x,y)&\!\mapsto\!&\left(\frac{(x^2-1)\,Q(x^{-2})^2}{x\,P(x^{-2})^2},\;
\frac{y\,Q(x^{-2})\,R(x^{-2})}{P(x^{-2})^3}\right),\\ \label{el4isog3}
(x,y)&\!\mapsto\!&\left(\frac{x^3\,Q(x^{-2})^2}{(x^{2}-1)\,P(x^{-2})^2},\;
\frac{x^5\,Q(x^{-2})\,R(x^{-2})}{y\,(x^2-1)\,P(x^{-2})^3}\right).
\end{eqnarray}
Indeed, a direct check that the transformed coordinate functions in (\ref{el4isog1})
satisfy the equation $y^2=x^3-x$ gives expression (\ref{pqrcheck1})
after the substitution $z\mapsto x^{-2}$. 
Similarly, transformations (\ref{el4isog2}) and (\ref{el4isog3}) are compatible with
the equation $y^2=x^3-x$ and, respectively, (\ref{pqr4check2}) or (\ref{pqr4check3}).
These endomorphisms fix the point at infinity, so they are isogenies. \qed\\

\noindent 
Let us refer to the degree $4n+1$, $4n+2$, $4n$ cases as, respectively,
degree $[4n+1]$, $[4n+2]$, $[4n]$ coverings.  
The addition law on $E_1$ can be translated to relations between the polynomials 
$P(z)$, $Q(z)$, $R(z)$ for different pull-back coverings. To reduce the number of cases,
we can represent a degree $[4n]$ triple $(P,Q,R)$ as a degree $[4n+2]$ type 
\mbox{$((1-z)P,Q,(1-z)R)$}; check this transformation between equations (\ref{pqr4check2})
and (\ref{pqr4check3}). Here is addition of two degree $[4n+1]$ coverings 
$(P_1,Q_1,R_1)$ and $(P_2,Q_2,R_2)$, represented as a degree $[4n+2]$ covering
$\left({\bf P}, {\bf Q}, {\bf R}\right)$:
\begin{eqnarray} \label{e4add}
\left({\bf P}, {\bf Q}, {\bf R}\right) &\equal &
\big( P_1^2Q_2^2-P_2^2Q_1^2,\;  P_1Q_1R_2-P_2Q_2R_1,\nonumber \\ &&
\;(1-z)(P_1^2Q_2^2+P_2^2Q_1^2)R_1R_2
-2P_1P_2Q_1Q_2(Q_1^2Q_2^2-zP_1^2P_2^2)\big).
\end{eqnarray}
The resulting polynomials ${\bf P}, {\bf Q}, {\bf R}$ may have common factors, and the triple may
simplify in a projective manner as $\left({\bf P}/h, {\bf Q}/h, {\bf R}/h^2\right)$. In fact, the polynomial
triple must typically have large common factors, as degree of the corresponding isogenies does not grow so fast. If ${\bf P}$ and ${\bf R}$ (but not ${\bf Q}$) are divisible by $(1-z)$, 
the "sum" actually has degree $[4n]$. Formula (\ref{e4add}) is void as a duplication formula
for $(P_2,Q_2,R_2)=(P_1,Q_1,R_1)$.


Similarly, here are formulas for "adding" two degree $[4n+2]$ coverings 
(represented as a degree $[4n]$ covering), and "adding" a degree $[4n+2]$ and a degree $[4n+1]$
covering (represented as a degree $[4n+1]$ covering), respectively:
\begin{eqnarray}  \label{e4add2}
\left({\bf P}, {\bf Q}, {\bf R}\right) &\!\equal\! &
\big( P_1^2Q_2^2-P_2^2Q_1^2,\;  P_1Q_1R_2-P_2Q_2R_1,\nonumber\\ &&
\;(P_1^2Q_2^2+P_2^2Q_1^2)R_1R_2
-2P_1P_2Q_1Q_2((1-z)^2Q_1^2Q_2^2-zP_1^2P_2^2)\big),\\ \label{e4add3}
\left({\bf P}, {\bf Q}, {\bf R}\right) &\!\equal\! &
\big( P_1^2Q_2^2-(1-z)P_2^2Q_1^2,\;  (1-z)P_1Q_1R_2-P_2Q_2R_1,\nonumber\\ &&
((1\!-\!z)P_1^2Q_2^2\!+\!P_2^2Q_1^2)R_1R_2
-2P_1P_2Q_1Q_2((1\!-\!z)Q_1^2Q_2^2\!-\!zP_1^2P_2^2)\big).
\end{eqnarray}
Expressions differ by placement of $(1-z)$ factors. Formally, a degree $[4n+2]$ case 
$(P,Q,R)$ can be represented in the degree $[4n+1]$ form
$\left( \sqrt{1-z}\,P,(1-z)Q,\sqrt{1-z}\,R\right)$.
If $z$ is considered fixed, we can recognize an addition law on the curve $(1-z)y^2=x^4-z$.

The action of the isogeny $\phi_i$ can be represented as $(P,Q,R)\mapsto (iP,Q,R)$,
which does not change the pull-back covering, as mentioned. However, this transformation
can be used to compute, say, a duplication formula via the action of $1+i$. 
In particular, the action of $1+i$ on degree $[4n+1]$ coverings
follows from formulas (\ref{e4add})--(\ref{e4add2}) with $(P_2,Q_2,R_2)=(iP_1,Q_1,R_1)$;
after the mentioned projective division by $(1-i)$ we get degree $[4n+2]$ and $[4n]$ coverings
\begin{eqnarray}
\left({\bf P}, {\bf Q}, {\bf R}\right) &\equal& \big( (1+i)P_1Q_1,\,R_1,\,Q_1^4+zP_1^4 \big),\\
\left({\bf P}, {\bf Q}, {\bf R}\right) &\equal& \big( (1+i)P_1Q_1,\,R_1,\,(1-z)^2Q_1^4+zP_1^4 \big).
\end{eqnarray}
Applying now the conjugate action of $1-i$ we get the duplication formulas, respectively,
\begin{eqnarray*}
\left({\bf P}, {\bf Q}, {\bf R}\right) &\equal& \big( 2P_1Q_1R_1,\,Q_1^4+zP_1^4,\,Q_1^8-6zP_1^4Q_1^4+z^2P_1^8 \big),\\
\left({\bf P}, {\bf Q}, {\bf R}\right) &\equal& \big( 2P_1Q_1R_1,\,(1-z)^2Q_1^4+zP_1^4,
(1-z)^4Q_1^8-6z(1-z)^2P_1^4Q_1^4+z^2P_1^8 \big).
\end{eqnarray*}
Here both resulting coverings have degree $[4n]$. By the isomorphism with isogenies,
the triples $(P,Q,R)$ defining our  transformation coverings can be given 
(modulo the mentioned projective equivalence)  a ring structure isomorphic to $\ZZ[i]$.
By induction on the norm of $u\in\ZZ[i]$ we observe
that the transformation corresponding to $u$ can be represented by a polynomial triple
evaluating at $z=0$ to $(u,1,1)$. Formulas (\ref{e4add}) and (\ref{e4add3}) with
$(P_2,Q_2,R_2)=(1,1,1)$ or $(i,1,1)$ 
give recursion relations for computing new polynomial triples with added 
$1$ or $i$ to the indexing Gaussian integer. We get the following polynomial triples
corresponding to the Gaussian integers $2+i$, $2+2i$, $3$:
\begin{eqnarray*} 
\left( 2+i-iz, 1+(2i-1)z, 1+(2-8i)z+z^2 \right),\\ 
\left( (2+2i)(1+z), 1-6z+z^2,1+20z-26z^2+20z^3+z^4\right),\\
\left( 3-6z-z^2, 1+6z-3z^2, 1-28z+6z^2-28z^3+z^4\right).
\end{eqnarray*}

\section{Elliptic integrals on the curve $y^2=x^3-1$}
\label{elints6}

Let $E_2$ denote the genus 1 curve defined by the equation $y^2=x^3-1$. We consider it as 
an elliptic curve in the usual way. 
Let $f_6=t^{-5/6}\,(1-t)^{-1/2}$ denote the integrand function of (\ref{elin6}). 
The algebraic curve with the function field $\CC(t,f_6)$ is isomorphic to $E_2$ by the isomorphisms
$(t,f_6)=\left(x^{-3},\,x^4y^{-1}\right)$ and $x=f_6^4\,t^3\,(1-t)^2$.
Therefore (\ref{elin6}) is an elliptic integral on $E_2$. 
By substituting $t=x^{-3}$ in the integral in (\ref{elin6}) we get:
\begin{equation} \label{elin6a}
\hpg{2}{1}{1/2,\,1/6}{7/6}{z}=\frac{z^{-1/6}}{2}
\,\int_{1/\!\sqrt[3]{z}}^{\infty}\;\frac{dx}{\sqrt{x^3-1}}.
\end{equation}
Like in the previous case, we show that transformations of the hypergeometric function in
(\ref{elin6a}) correspond to isogeny endomorphisms of $E_2$.

The elliptic curve $E_2$ has complex multiplication by a cubic roots of unity.
Let $\omega$ is a primitive cubic root of unity; then the ring of isogenies of $E_2$
is isomorphic to $\ZZ[\omega]$.
We identify $\omega\in\ZZ[\omega]$ with the automorphism $(x,y)\mapsto(\omega\,x,y)$ on $E_2$.
For $a+b\omega\in\ZZ[\omega]$ let $\phi_{a+b\omega}=(\varphi_{a+b\omega},\psi_{a+b\omega})$
denote the corresponding isogeny; its degree equal to the norm
$a^2-ab+b^2$ of $a+b\omega$.

\begin{theorem} \label{th:el6tr}
Let $\phi:E_2\to E_2$ denote a non-zero isogeny on $E_2$, of degree $d$.
The isogeny transforms the $x$-coordinate as $x \mapsto x\,\mu(x^3)$,
where $\mu(s)$ is a rational function with a finite non-zero limit 
$\mu_0=\lim_{s\to\infty}\mu(s)$. We have the identity
\begin{equation} \label{el6tr}
\hpg{2}{1}{1/2,\,1/6}{7/6}{\,z} = \left(\frac{\mu(1/z)}{\mu_0}\right)^{-1/2}\;
\hpg{2}{1}{1/2,\,1/6}{7/6}{\frac{z}{\mu\!\left(1/z\right)^3}},
\end{equation}
The degree of the rational function $z/\mu\!\left(1/z\right)^3$ is equal to $d$.
\end{theorem}
\proof As in the proof of Theorem \ref{th:el4tr}, 
let $\varphi(x)=\varphi(x,y)$ and $\psi(x)=\psi(x,y)$ denote 
the $x,y$ components of the isogeny $\phi$. 

This integrand in (\ref{elin6a}) is a holomorphic differential $1$-form on
$E_2$. The substitution $x\mapsto \varphi(x,y)$ into the integral in (\ref{elin4a}) must be
an integral of  a holomorphic differential $1$-form on $E_2$ again. 
Since the linear space of holomorphic differentials is one-dimensional,
the transformed differential form must be proportional to $dx/\sqrt{x^3-1}$. 
The upper integration bound does not change; the lower
bound is transformed as $z\mapsto\varphi(z^{-1/3})^{-3}$.
Once we show that $\varphi(x)=x\,\mu(x^3)$ as stated in theorem, transformation (\ref{el6tr})
follows. 

The addition law on $E_2$ is defined by formula (\ref{elladd}) with $u,v\in\ZZ[\omega]$ and
\begin{equation}
\lambda=\frac{\psi_u-\psi_v}{\varphi_u-\varphi_v}\qquad\mbox{or}\qquad
\lambda=\frac{\varphi_u^2+\varphi_u\varphi_v+\varphi_v^2}{\psi_u+\psi_v}.
\end{equation}
By induction on the norm $|u|$ we conclude that any isogeny can be written
in the form $\left(x\mu(x^3),y\nu(x^3)\right)$, where $\mu(s)$, $\nu(s)$ are
rational functions satisfying $\lim_{s\to\infty} \mu(s)=1/u^2$ and
$\lim_{s\to\infty} \nu(s)=1/u^3$. Transformation (\ref{el6tr}) follows. 

For the last statement, note that the points $(x,y)=(0,\pm i)$ are in the kernel
of all isogenies $\phi_u$ with $d=|u|$ divisible by 3; hence $s=0$ is a pole of $\mu(s)$
when $d$ is divisible by 3. The degree of the rational function $\varphi(x)=x\,\mu(x^3)$ is $d$, 
and the degree of $\mu(s)$ is equal to $\lfloor d/3\rfloor$ for any $d$.
The factor $z$ in $z/\mu(1/z)^3$ simplifies in the same way, exactly when $d$ is divisible by 3.
\qed\\

\noindent
Here are the explicit transformations corresponding to the algebraic integers
$1-\omega$, $2$, $3$ and $3\omega+1$, respectively:
\begin{eqnarray*}
\hpg{2}{1}{\!1/2,\,1/6}{7/6}{z} & \equal & \frac{1}{\sqrt{1-4z}}\;
\hpg{2}{1}{\!1/2,\,1/6}{7/6}{\frac{27\,z}{(4z-1)^3}},\\
\hpg{2}{1}{\!1/2,\,1/6}{7/6}{z} & \equal & \sqrt{\frac{1-z}{1+8z}}\;
\hpg{2}{1}{\!1/2,\,1/6}{7/6}{\frac{64\,z\,(1-z)^3}{(8z+1)^3}},\\
\hpg{2}{1}{\!1/2,\,1/6}{7/6}{z} & \equal &
\frac{1-4z}{\sqrt{1\!+\!96z\!+\!48z^2\!-\!64z^3}}\,\hpg{2}{1}{\!1/2,1/6}{7/6}
{\frac{-729\,z\,(4z-1)^6}{(64z^3\!-\!48z^2\!-\!96z\!-\!1)^3}}.\\
\hpg{2}{1}{\!1/2,1/6}{7/6}{z} &\equal&\frac{1-4z/(3\omega\!+\!1)}
{\sqrt{1\!-\!(44\!+\!48\omega)z\!+\!(48\omega\!+\!16)z^2}}\times\nonumber\\
&&\;\hpg{2}{1}{\!1/2,1/6}{7/6}{\frac{z\;(4z-3\omega\!-\!1)^6}
{((48\omega\!+\!16)z^2\!-\!(44\!+\!48\omega)z\!+\!1)^3}}.
\end{eqnarray*}
The first two formulas are special cases of classical transformations.
The isogenies corresponding to the roots of unity in $\ZZ[\omega]$
give the trivial transformation of (\ref{elin6a}),  so the group of 
transformations (\ref{el6tr}) is isomorphic to the multiplicative
$\ZZ[\omega]^*/\left(\pm 1,\pm\omega,\pm\omega^{-1}\right)$.

\begin{theorem}
Any transformation of the form $(\ref{el6tr})$ comes from an isogeny endomorphism of $E_2$
as described in Theorem $\ref{th:el6tr}$.
\end{theorem}
\proof There are four entries in Table \ref{elltab} giving explicit ramification patterns of 
possible pull-back coverings of Euler's hypergeometric equation with the local exponent
differences $(1/2,1/3,1/6)$ into itself. For identity $(\ref{el6tr})$ we should have the local exponent
difference $1/6$ at the points $z=0$ on both projective lines, lying above each other.
The points $z=\infty$ have the local exponent difference $1/3$.

With this setting, pull-back coverings of degree $6n+1$ and $6n+3$  imply the polynomial identity
\begin{equation} \label{pqr6check1}
(1-z)\,R(z)^2=Q(z)^3-z\,P(z)^6,
\end{equation}
where polynomials $P(z),Q(z),R(z)$ have degree $n,2n,3n$ (for the covering degree $6n+1$) or
$n,2n+1,3n+1$ (for the covering degree $6n+3$). Similarly, coverings of degree $6n+4$ and $6n$
imply the polynomial identity
\begin{equation} \label{pqr6check2}
R(z)^2=Q(z)^3-z\,(1-z)^3\,P(z)^6,
\end{equation}
where polynomials $P(z),Q(z),R(z)$ have degree $n,2n+1,3n+2$ (for the covering degree $6n+4$) 
or $n-1,2n,3n$ (for the covering degree $6n$). 

From identities (\ref{pqr6check1}) and (\ref{pqr6check2}) we get the following 
endomorphisms of $E_2$, respectively:
\begin{eqnarray} \label{el6isog1}
(x,y)&\!\mapsto\!&\left(\frac{x\,Q(x^{-3})}{P(x^{-3})^2},\;
\frac{y\;R(x^{-3})}{P(x^{-3})^3}\right),\\ \label{el6isog2}
(x,y)&\!\mapsto\!&\left(\frac{x^4\,Q(x^{-3})}{(x^3-1)\,P(x^{-3})^2},\;
\frac{x^6\,R(x^{-3})}{y\,(x^3-1)\,P(x^{-3})^3}\right).
\end{eqnarray}
The endomorphisms fix the point at infinity, so they are isogenies. \qed\\

\noindent
Like in the previous section, the addition law on $E_2$ can be translated to relations between
the polynomials  $P(z)$, $Q(z)$, $R(z)$ for different pull-back coverings. To reduce the number of cases,
we can represent a triple $(P,Q,R)$ for even degree $6n+4$ or $6n$ coverings as an odd degree
$6n+1$ or $6n+3$ case $\big((1-z)P,(1-z)Q,(1-z)R\big)$. 
Here is addition of two coverings $(P_1,Q_1,R_1)$ and $(P_2,Q_2,R_2)$ of odd degree, 
represented as an odd degree covering:
\begin{eqnarray}  \label{e6add}
\hspace{-28pt}&&\big( P_1^2Q_2-P_2^2Q_1,\;  
Q_1Q_2(P_1^2Q_2+P_2^2Q_1)-2(1\!-\!z)P_1P_2R_1R_2-2zP_1^4P_2^4,\nonumber\\ 
\hspace{-28pt}&&\;P_1Q_2^2R_1(P_1^2Q_2\!+\!3P_2^2Q_1)-
P_2Q_1^2R_2(P_2^2Q_1\!+\!3P_1^2Q_2)
+4zP_1^3P_2^3(P_1^3R_2\!-\!P_2^3R_1)\big).
\end{eqnarray}
The resulting polynomial triple may
simplify in a projective manner as $\left({\bf P}/h, {\bf Q}/h^2, {\bf R}/h^3\right)$.
By the isomorphism with isogenies, the polynomial triples can be given a ring structure isomorphic 
to $\ZZ[\omega]$. By induction on the norm of $u\in\ZZ[\omega]$ we observe that
the pull-back transformation corresponding to $u$ 
can be represented by the a polynomial triple evaluating at $z=0$ to $(u,1,1)$.
Multiplication by the roots of unity in $\ZZ[\omega]$ is just as simple as in the 
$(1/2,1/4,1/4)$ case. In particular, multiplication of $(P,Q,R)$ by $1+\omega$ gives
$\left( (1+\omega)P, Q, R \right)$, as can be checked using (\ref{e6add}) with
$(P_2,Q_2,R_2)=(\omega P_1,Q_1,R_1)$. Multiplication of odd degree $(P,Q,R)$ 
by $1-\omega$ gives
$\left( (1-\omega)PQ, Q^3-4zP^6, (Q^3+8zP^6)R \right)$. Since $2=(1+\omega)+(1-\omega)$,
the duplication formula for odd degree $(P,Q,R)$ is
\begin{equation}
\left( 2PR,\; Q(Q^3+8zP^6),\; Q^6-20zP^6Q^3-8z^2P^{12} \right),
\end{equation}
represented as an even degree case. 
Formula (\ref{e6add}) with $(P_2,Q_2,R_2)=(1,1,1)$ or $(\omega,1,1)$ 
gives recursion relations, on the lattice $\ZZ[\omega]$, between the polynomial triples.
One may also use the conjugation $\omega\mapsto-\omega-1$. 
Here are the triples corresponding to $2-\omega$ and $2-2\omega$, respectively:
\begin{eqnarray*}
\big( 2-\omega+(4+4\omega)z, 1-(44+48\omega)z+(16+48\omega)z^2, \hspace{42pt} \\
1+(96+108\omega)z+(48-432\omega)z^2-64z^3 \big),\\
\big((2-2\omega)(1+8z), (1-4z)(1-228z+48z^2-64z^3), \hspace{42pt} \\
(1-20z-8z^2)(1+536z-1344z^2+2048z^3-512z^4)\big).
\end{eqnarray*}

As Table \ref{elltab} indicates, there are transformations of hypergeometric functions
with the local exponent differences $(1/2,1/3,1/6)$ to hypergeometric equations with the
local exponent differences $(1/3,1/3,1/3)$ and $(1/6,1/6,2/3)$. We consider these transformations
in Section \ref{otherells}.

\section{Elliptic integrals on the curve $X^3+Y^3=1$}
\label{elints3}

Let $E_3$ denote the genus 1 curve defined by the equation
$X^3+Y^3=1$. It is is isomorphic to $E_2$ via the isomorphisms
\begin{equation} \label{e2e3iso}
(X,\,Y)\mapsto \left(\frac{2^{2/3}}{X+Y},\;\sqrt{3}\,\frac{X-Y}{X+Y}\right),
\qquad (x,\,y)\mapsto\left( \frac{y+\sqrt{3}}{2^{1/3}\sqrt{3}\,x},\;
\frac{\sqrt{3}-y}{2^{1/3}\sqrt{3}\,x}\right).
\end{equation}
Under this ismorphism, $E_3$ can be considered as an elliptic curve 
with the point $(X,Y,1)=(1,-1,0)$  as the neutral element of its additive group. 
Then the additive opposite of $(X,Y)\in E_3$ is the point $(Y,X)$,
and the complex multiplication by $\omega\in\ZZ[\omega]$ is
the isogeny $(X,Y)\mapsto(\omega^{-1}X,\omega^{-1}Y)$.

Let $f_3=t^{-2/3}\,(1-t)^{-2/3}$ be the integrand of (\ref{elin3}),
and let $G$ be the algebraic curve with the function field
$\CC(t,f_3)$. An isomorphism between $E_2$ and $G$ can be given by 
\begin{eqnarray}
(t,f_3)\mapsto\left( 2^{2/3}f_3\,t(1-t), i(2t-1) \right),\qquad
(x,y)\mapsto\left( \frac{1-iy}2, \, \frac{2^{4/3}}{x^2} \right).
\end{eqnarray}
The compatible isomorphism between $E_3$ and $G$ is given by 
\begin{eqnarray*}
\;(X,Y) \mapsto\! \left( \frac{(\omega\!+\!1)Y-\omega X}{X+Y},\, (X+Y)^2 \right),\ 
(t,f_3) \mapsto\! \left( \frac{1+\omega-t}{\sqrt{-3}\,f_3t(t\!-\!1)},\,
\frac{\omega+t}{\sqrt{-3}\,f_3t(t\!-\!1)} \right).
\end{eqnarray*}
Here $\sqrt{-3}=2\omega+1$. 

We see that (\ref{elin3}) is an elliptic integral. A convenient substitution into the integral is 
$t=X^{-3}$. We get the following integral of a holomorphic form on $E_3$: 
\begin{equation} \label{elin3a}
\hpg{2}{1}{1/3,\,2/3}{4/3}{z} = z^{-1/3}\,
\int_{1/\!\sqrt[3]{z}}^\infty\;\frac{dX}{(X^3-1)^{2/3}}.
\end{equation}
The morphism $E_3\to G$ given by $(X,Y)\mapsto
\left(X^{-3},\,X^4\,Y^{-2}\right)$ has degree 3, as can be checked using Gr\"oebner bases techniques.
This morphism factors as a composition of the following morphism $E_3\to E_2$ and
isomorphism $E_2\to G$: 
\begin{equation} \label{e2e3m3}
(X,Y)\mapsto\left(2^{2/3}XY,\,i(2X^3-1)\right),\qquad
(x,y)\mapsto \left( \frac2{1-iy}, \frac{(1-iy)^2}{2^{2/3}x^2} \right).
\end{equation}
The degree 3 morphism $E_3\to E_2$ sends all points of $E_3$ at infinity
to the origin of $E_2$.

We show that transformations of the hypergeometric function in
(\ref{elin3}) correspond to the isogenies of $E_3$, similarly to
the previous cases. 
\begin{theorem} \label{th:el3tr}
Let $\phi:E_3\to E_3$ denote a non-zero isogeny on $E_3$, of degree $d$.
The isogeny either transforms the $X$- or $Y$-coordinate to $X\mu(X^3)$,
or transforms either coordinate to $\eta(X^3)/XY$,
where $\mu(s)$ or $\eta(s)$ is a rational function with a finite non-zero limit 
$\mu_0=\lim_{s\to\infty}\mu(s)$ or $\eta_0=\lim_{s\to\infty}\eta(s)$. 
We have the identity
\begin{equation} \label{el3tra}
\hpg{2}{1}{1/3,\,2/3}{4/3}{\,z} = \frac{\mu_0}{\mu(1/z)}\;
\hpg{2}{1}{1/3,\,2/3}{4/3}{\frac{z}{\mu\!\left(1/z\right)^3}}
\end{equation}
or the identity
\begin{equation} \label{el3trb}
\hpg{2}{1}{1/3,\,2/3}{4/3}{\,z} = \frac{(1-z)^{1/3}\,\eta_0}{\eta(1/z)}\;
\hpg{2}{1}{1/3,\,2/3}{4/3}{\frac{z\,(1-z)}{\eta\!\left(1/z\right)^3}}.
\end{equation}
The degree of the rational function $z/\mu\!\left(1/z\right)^3$ or $z(1\!-\!z)/\eta\!\left(1/z\right)^3$
is equal to $d$.
\end{theorem}
\proof Like in the proof of Theorem \ref{th:el6tr}, 
let $\varphi(X)=\varphi(X,Y)$ and $\psi(x)=\psi(X,Y)$ denote 
the $X,Y$ components of the isogeny $\phi$. 

This integrand in (\ref{elin3a}) is a holomorphic differential $1$-form on
$E_3$. The substitution $X\mapsto \varphi(X,Y)$ into the integral in (\ref{elin4a}) must be
an integral of  a holomorphic differential $1$-form on $E_3$ again. 
Since the linear space of holomorphic differentials is one-dimensional,
the transformed differential form must be proportional to $(X^3-1)^{-2/3}dX$. 
The upper integration bound does not change; the lower
bound is transformed as $z\mapsto\varphi(z^{-1/3})^{-3}$.
Here is an explicit expression for 
addition of two isogenies $(\varphi_{u},\psi_{u})$ and $(\varphi_{v},\psi_{v})$ on $E_3$:
\begin{equation} \label{add3uv}
(\varphi_{u+v},\psi_{u+v})=\left(
\frac{\varphi_u\psi_v^2-\varphi_v\psi_u^2}{\varphi_u\psi_u-\varphi_v\psi_v},\,
\frac{\varphi_v^2\psi_u-\varphi_u^2\psi_v}{\varphi_u\psi_u-\varphi_v\psi_v}\right).
\end{equation}
Using induction on $|u|$ one shows that any isogeny $(\varphi_u,\psi_u)$ of $E_3$ has one of the
following three forms: 
\begin{itemize}
\item $\left(X\mu(X^3),Y\nu(X^3)\right)$, where $\lim_{s\to\infty} \mu(s)=\lim_{s\to\infty} \nu(s)=1/u$.
\item $\left(Y\nu(X^3),X\mu(X^3)\right)$, where $\lim_{s\to\infty} \mu(s)=\lim_{s\to\infty} \nu(s)=-1/u$.
\item $\left(\mu(X^3)Y^2/X,\nu(X^3)X^2/Y\right)$, where 
$\lim_{s\to\infty} \mu(s)=\lim_{s\to\infty} \nu(s)=1/u$.
\end{itemize}
The form depends on the residue of $u$ modulo the lattice generated by $3$ and $\omega+2$.
Regarding the last case, it is useful to observe that 
\[
\frac{Y^2}{X}=\frac{1-X^3}{XY}=\frac{1-X^3}{X^3}\cdot \frac{X^2}{Y}.
\] 
Transformations (\ref{el3tra})--(\ref{el3trb}) and their degree easily follow.
\qed\\

\noindent
Here are a few explicit examples transformations that come from
isogenies of$E_3$. The corresponding isogenies
$(\varphi_u,\psi_u)$ have $u=1-\omega$, $3$ or $3+\omega$, as in
the previous section. The first identity is a special case of a classical
(though not well-known) cubic transformation. 
\begin{eqnarray*}
\hpg{2}{1}{1/3,2/3}{4/3}{z} & \equal &
\frac{(1-z)^{1/3}}{1+\omega^2z}\;\hpg{2}{1}{1/3,\,2/3}{4/3}
{\frac{3\,(2w+1)\,z\,(z-1)}{(z+\omega)^3}}.\\
\hpg{2}{1}{1/3,2/3}{4/3}{z} & \equal &
\frac{(1\!-\!z\!+\!z^2)\,(1\!-\!z)^{1/3}}
{1\!+\!3z\!-\!6z^2\!+\!z^3}\;\hpg{2}{1}{1/3,\,2/3}{4/3}
{\frac{27\,z\,(z-1)\,(z^2-z+1)^3}{(z^3-6z^2+3z+1)^3}}.\\
\hpg{2}{1}{\!1/3,2/3}{4/3}{z}&\equal&
\frac{1-z-z^2/(3\omega\!+\!2)}{1\!+\!(3\omega\!+\!2)z\!-\!(3\omega\!+\!2)z^2}\\
&&\times \hpg{2}{1}{\!1/3,2/3}{4/3}{\frac{z\;(z^2+(3\omega\!+\!2)z-3\omega\!-\!2)^3}
{(1+(3\omega\!+\!2)z-(3\omega\!+\!2)z^2)^3}}.
\end{eqnarray*}

\begin{theorem}
Any transformation of the form $(\ref{el3tra})$ or $(\ref{el3trb})$ comes from an isogeny
endomorphism of $E_3$ as described in Theorem $\ref{th:el3tr}$.
\end{theorem}
\proof There are two entries in Table \ref{elltab} giving explicit ramification patterns of 
possible pull-back coverings of Euler's hypergeometric equation with the local exponent
differences $(1/3,1/3,1/3)$ into itself. For a hypergeometric identity we should specify the 
points $z=0$ on both projective lines to lie above each other.

With this setting, pull-back coverings of degree $3n+1$ imply the polynomial identity
\begin{equation} \label{pqr3check1}
z\,P(z)^3+(1-z)\,R(z)^3=Q(z)^3,
\end{equation}
where polynomials $P(z),Q(z),R(z)$ have degree $n$. Pull-back coverings of degree $3n$
imply the polynomial identity
\begin{equation} \label{pqr3check2}
z\,(z-1)\,P(z)^3=Q(z)^3+R(z)^3,
\end{equation}
where polynomials $P(z),Q(z),R(z)$ have degree $n-1,n,n$, respectively.
We get the following endomorphisms of $E_3$, respectively:
\begin{eqnarray} \label{el3isog1}
(X,Y)&\!\mapsto\!&\left(\frac{X\,Q(X^{-3})}{P(X^{-3})},\;
\frac{Y\,R(X^{-3})}{P(X^{-3})}\right),\\ \label{el3isog2}
(X,Y)&\!\mapsto\!&\left(\frac{X^2\,Q(X^{-3})}{Y\,P(X^{-3})},\;
\frac{X^2\,R(X^{-3})}{Y\,P(X^{-3})}\right).
\end{eqnarray}
The endomorphisms fix the point at infinity, so they are isogenies. \qed\\

\noindent
Like in the previous sections, the addition law on $E_3$ can be translated to relations between
the polynomials  $P(z)$, $Q(z)$, $R(z)$ for different pull-back coverings. The form of these identities
follows from the addition formula (\ref{add3uv}) straightforwardly.

\section{Other transformations of elliptic integrals}
\label{otherells}

Elliptic integrals (\ref{elin6}) and (\ref{elin3}) are defined on
isomorphic curves $E_2$ and $E_3$. Therefore we expect that these
two integrals are related. Indeed, we have the following classical
quadratic transformation:
\begin{equation} \label{e3e6tr}
\hpg{2}{1}{1/3,\,2/3}{4/3}{\,z}=(1-z)^{-1/6}\;\hpg{2}{1}{1/2,\,1/6}
{7/6}{\frac{z^2}{4(z-1)}}.
\end{equation}
This formula can also be derived by substituting $x\mapsto
2^{2/3}X(1-X^3)^{1/3}$ in (\ref{elin6a}). The corresponding
morphism $E_3\to E_2$ can be given by
\begin{equation} \label{e3e6isog}
(X,Y)\mapsto (x,y)= \left(2^{2/3}XY,\,i-2iX^3\right).
\end{equation}
Notice that this morphism has degree 3 whereas transformation
(\ref{e3e6tr}) is quadratic.

\begin{theorem}
Any transformation of a hypergeometric equation with the local exponent differences
$(1/2,1/3,1/6)$ to a hypergeometric equation with the local exponent differences
$(1/3,1/3,1/3)$ are compositions of quadratic transformation $(\ref{e3e6tr})$
and an isogeny endomorphism of $E_2$ or $E_3$ described in Sections
$\ref{elints6}$ and $\ref{elints3}$.
\end{theorem}
\proof We reiterate that algebraic transformations 
of hypergeometric functions
transform their hypergeometric equations as well. From Table \ref{elltab}
it is straightforward to conclude that a
transformation between hypergeometric equations for (\ref{elin6})
and (\ref{elin3}) induce an isogeny $\phi:E_3\to E_2$ that maps
the three infinite points of $E_3$ to the infinite point of $E_2$.
We have to check every possibility from Table \ref{elltab}. To fix
the ideas, consider a morphism of degree $6n+2$. Existence of such
a morphism implies the identity
\begin{equation} \label{pqrcheck2}
R_{3n+1}(z)^2=(1-z)Q_{2n}(z)^3-z^2\,P_n(z)^6,
\end{equation}
where $P_n(z)$, $Q_{2n}(z)$ and $R_{3n+1}(z)$ are some polynomials
of degree $n$, $2n$ and $3n+1$ respectively. Then the following
transformation is a morphism from $E_3$ to $E_2$:
\begin{equation}
(X,Y)\mapsto (x,y)=\left(
\frac{X\,Y\,Q_{2n}(X^{-3})}{P_n(X^{-3})^2},\;
\frac{X^3\;R_{3n+1}(X^{-3})}{P_n(X^{-3})^3} \right).
\end{equation}
Moreover, this is an isogeny since the neutral element of $E_3$ is
mapped to the neutral element of $E_2$. Also the infinite points
of $E_3$ are mapped to the infinite point of $E_2$. The same claim
can be checked for other possibilities in a similar way.

Now we identify the isogenies on $E_3$ and $E_2$ using isomorphism
(\ref{e2e3iso}). Recall that these isogenies form the ring
isomorphic to $\ZZ[\omega]$. In particular, isogeny
(\ref{e3e6isog}) is then identified with $\sqrt{-3}$, which is
$\pm(1+2\omega)$. Consider an isogeny $\phi$ that is induced by a
hypergeometric transformation. Since the infinite points of $E_3$
have order 3 (or 1) in the additive group, and since $(1+2\omega)$
is the only prime ideal above $(3)$ in $\ZZ[\omega]$, the isogeny
$\varphi$ can be factored as $\varphi_1\circ\varphi_2$, where
$\varphi_1$ is equal to (\ref{e3e6isog}). By considering
$\varphi_2$ as an isogeny on $E_2$ we can express the
corresponding hypergeometric identity as a composition of
(\ref{e3e6tr}) and a transformation from Subsection \ref{elints6}.
(Other factorizations and compositions are possible.) 
\qed\\

As a corollary, transformations of
degree $6n+4$ between (\ref{elin6}) and (\ref{elin3}) suggested by
Table \ref{elltab} are non-existent.


It remains to consider transformations of the hypergeometric
equation for (\ref{elin6}) to the hypergeometric equation with the
local exponent differences $(2/3,1/6,1/6)$. A solution of the
latter equation is the following:
\begin{equation} \label{hpelin}
\hpg{2}{1}{1/3,\,1/6}{7/6}{z} = \frac{z^{-1/6}}{6}\, \int_0^z
t^{-5/6}\,(1-t)^{-1/3}\,dt.
\end{equation}
Let $f=t^{-5/6}\,(1-t)^{-1/3}$ be the integrand of (\ref{hpelin}),
and let $H$ be the algebraic curve with the function field
$\CC(t,f)$. This curve is isomorphic to the hyperelliptic curve
$Y^2=X^6+1$ by the isomorphisms $t=\left(Y-X^3\right)^2$ and
\begin{equation}
(X,Y)=\left( \frac{1}{2^{1/3}\,tf},\;
\frac{1}{2\,t^3f^3}\frac{1+t}{1-t} \right).
\end{equation}
Therefore (\ref{hpelin}) is a hyperelliptic integral. The
substitution $t\to \left(\sqrt{X^6+1}-X^3\right)^2$ transforms
this integral into
\begin{equation} \label{hpelin2}
\hpg{2}{1}{1/3,\,1/6}{7/6}{\,z} = \frac{z^{-1/6}}{2^{1/3}}\,
\int_{\theta(z)}^{\infty} \frac{X\,dX}{\sqrt{X^6+1}},
\qquad\mbox{where}\quad
\theta(z)=\frac{(1-z)^{1/3}}{2^{1/3}\,z^{1/6}}.
\end{equation}
If $\varphi(X)$ is the $x$-component of an algebraic map $H\mapsto
E_2$, then there is a transformation of the hypergeometric
function with the argument $\varphi(\theta(z))^{-3}$. The simplest
transformation between (\ref{elin6}) and (\ref{hpelin}) is the
following:
\begin{equation} \label{hpelltr}
\hpg{2}{1}{1/3,\,1/6}{7/6}{z}=(1-z)^{-1/3}\;\hpg{2}{1}{1/2,\,1/6}
{7/6}{-\frac{4z}{(z-1)^2}}.
\end{equation}
The corresponding morphism $H\to E_2$ can be given by
$(X,Y)\mapsto \left(-X^2,\,iY\right)$.
\begin{theorem}
Any transformation of a hypergeometric equation with the local exponent differences
$(1/2,1/3,1/6)$ to a hypergeometric equation with the local exponent differences
$(1/3,1/6,1/6)$ are compositions of quadratic transformation $(\ref{hpelltr})$
and an isogeny endomorphism of $E_2$ described in Section $\ref{elints6}$.
\end{theorem}
\proof From Table \ref{elltab} we conclude that such a pull-back covering 
gives a pull-back covering $H\to E_2$, that moreover maps the (singular) infinite point of $H$ to the
infinite point of $E_2$. The covering must factor via $(X,Y)\mapsto\left(-X^2,\,iY\right)$
because the differential $XdX/\sqrt{X^6+1}$ is invariant under the involution $X\mapsto -X$.
\qed

\section{Connection formulas}
\label{appendix}

The hypergeometric functions in (\ref{elin4})--(\ref{elin3}) fall into the category of
degenerate Gauss hypergeometric functions of \cite[Section 4]{degeneratehpg} with $n=0$.
In general, Euler's differential equation for hypergeometric functions in this category
has terminating solutions, non-abelian monodromy group and no logarithmic points.
In the special case $n=0$, twelve of the 24 Kummer's solutions represent trivial 
(i.e., constant or power) functions.

In general, Euler's differential equation for $\displaystyle\hpg21{-n,\,c\,}{1-a}{\,z\,}$ with $n=0$ is
\begin{equation} \label{eq:ehpg}
z\,(1-z)\,\frac{d^2y(z)}{dz^2}+\big(1-a-(1+c)\,z\big)\,\frac{dy(z)}{dz}=0.
\end{equation}
We assume $a, b\in\CC\setminus\ZZ$. The local exponents at $z=0$ are $0,a$; 
at $z=\infty$ are $0,c$; and at $z=1$ are $0,1-a-c$. 
A general solution of (\ref{eq:ehpg}) can be written as follows:
\begin{equation}
C_1\int z^{a-1}\,(1-z)^{-a-c}\,dz+C_2
\end{equation}
Let us define $b=1-a-c$. The 12 non-trivial Kummer's series solutions of (\ref{eq:ehpg})
represent the following three different functions:
\begin{eqnarray*} \label{eu4a}
z^a\;\hpg{2}{1}{a,\,1-b\,}{1+a}{\,z} &=& z^a\,(1-z)^{b}\;\hpg{2}{1}{1,\,a+b}{1+a}{\,z}\\
&=&z^a\,(1-z)^{b-1}\;\hpg{2}{1}{1,\,1-b\,}{1+a}{\,\frac{z}{z-1}}\\ \label{eu4c}
&=&z^a\,(1-z)^{-a}\;\hpg{2}{1}{a,\,a+b}{1+a}{\,\frac{z}{z-1}}, 
\\ 
(z-1)^{b}\;\hpg{2}{1}{b,\,1-a}{1+b}{1-z} \! &=& z^a\,(z-1)^{b}\;\hpg21{1,\,a+b}{1+b}{1-z} \\
&=&z^{a-1}\,(z-1)^{b}\;\;\hpg{2}{1}{1,\,1-a\,}{1+b}{1-\frac1{z}} \\
&=&z^{-b}\,(z-1)^{b}\;\hpg{2}{1}{b,\,a+b}{1+b}{1-\frac1{z}}, 
\\ 
z^{-c}\;\hpg21{c,\,1-b}{1+c}{\,\frac1z\,}  &=& z^{a-1}\,(z-1)^{b}\;\hpg{2}{1}{1,\,b+c}{1+c}{\,\frac1z\,}\\
&=&z^a\,(z-1)^{b-1}\;\hpg{2}{1}{1,\,1-b\,}{1+c}{\,\frac{1}{1-z}}\\ \label{eu4c}
&=&(z-1)^{-c}\;\hpg{2}{1}{c,\,b+c}{1+c}{\,\frac{1}{1-z}}.
\end{eqnarray*}
If $\mbox{Re}(a)>0$, $\mbox{Re}(b)>0$, $\mbox{Re}(c)>0$, 
the three functions can be written as the following integrals, respectively:
\begin{eqnarray}
a\int_0^z t^{a-1}\,(1-t)^{b-1}\,dx, \quad b\int_1^z t^{a-1}\,(t-1)^{b-1}\,dx, \quad
c\int_z^\infty t^{a-1}\,(t-1)^{b-1}\,dx.
\end{eqnarray}
Here are connection formulas for the hypergeometric solutions:
\begin{eqnarray} \label{eq:con01}
\frac{z^a}{a}\;\hpg{2}{1}{a,\,1-b}{1+a}{\,z}+\frac{(1-z)^{b}}{b}\;\hpg{2}{1}{b,\,1-a}{1+b}{1-z}=
\frac{\Gamma(a)\Gamma(b)}{\Gamma(a+b)},\\
\frac{(-z)^{a}}{a}\;\hpg{2}{1}{a,\,1-b}{1+a}{\,z}+\frac{(-z)^{-c}}{c}\;\hpg{2}{1}{c,\,1-b}{1+c}{\,\frac1z\,}=
\frac{\Gamma(a)\Gamma(c)}{\Gamma(a+c)},\\
\frac{(z-1)^{b}}{b}\;\hpg{2}{1}{b,\,1-a}{1+b}{1-z}+\frac{z^{c}}c\;\hpg{2}{1}{c,\,1-b}{1+c}{\,\frac1z\,}=
\frac{\Gamma(b)\Gamma(c)}{\Gamma(b+c)}.
\end{eqnarray}
The first formula is valid on the complex plane cut along  $(1,\infty)$ and $(-\infty,0)$;
the second formula is valid on $\CC\setminus(0,\infty)$; 
and the last formula is valid on $\CC\setminus(-\infty,1)$. 

Let us denote $\displaystyle F=\frac{z^a}{a}\;\hpg{2}{1}{a,\,1-b}{1+a}{\,z}$. 
Analytic continuation of $F$ along paths around the singularities $z=0$, $z=1$
determines the {\em monodromy group}. In general, analytic continuation to 
the same point changes $F$ to a solution $C_1F+C_2$ of differential equation (\ref{eq:ehpg}).
Specifically, the action on $F$ is the following:
\begin{equation}
\sigma_0F=\exppi{a}\,F, \qquad 
\sigma_1F=\exppi{b}\,F+\left(1-\exppi{b}\right)\frac{\Gamma(a)\Gamma(b)}{\Gamma(a+b)},
\end{equation}
where $\sigma_0$ represents analytic continuation along an anti-clockwise path around $z=0$ (only),
and $\sigma_1$ represents analytic continuation along clockwise path around $z=1$. 
In particular, analytic continuation along the Pochhammer path 
$\sigma^{-1}_0\sigma^{-1}_1\sigma_0\sigma_1$ gives 
\begin{equation} \label{eq:pochha}
\sigma^{-1}_0\sigma^{-1}_1\sigma_0\sigma_1F=F+
\left(1-\exppi{a}\right)\left(1-\exppi{b}\right)\,\frac{\Gamma(a)\Gamma(b)}{\Gamma(a+b)}.
\end{equation}

In the case $(a,b,c)=(1/4,1/2,1/4)$ with $z=x^2$, connection formula (\ref{eq:con01}) becomes
\begin{eqnarray} 
2\sqrt{x}\;\hpg{2}{1}{1/2,\,1/4}{5/4}{\,x^2}+\sqrt{1-x^2}\;\hpg{2}{1}{1/2,\,3/4}{3/2}{1-x^2}=
\frac{\Gamma(1/4)^2}{2\,\sqrt{2\pi}}.
\end{eqnarray}
However, it is valid only on the right half $\mbox{Re}(x)>0$ of the complex $x$-plane cut along 
$(1,\infty)$, as the imaginary axis is the pre-image of the $z$-cut $(0,-\infty)$ with respect to $z=x^2$.
We have to use other branch of $z^{1/4}$ in (\ref{eq:con01}) on the left half $\mbox{Re}(x)<0$ of the complex $x$-plane. The following formula holds on the left half-plane cut along $(-1,-\infty)$:
\begin{eqnarray} 
2i\sqrt{x}\;\hpg{2}{1}{1/2,\,1/4}{5/4}{\,x^2}+\sqrt{1-x^2}\;\hpg{2}{1}{1/2,\,3/4}{3/2}{1-x^2}=
\frac{\Gamma(1/4)^2}{2\,\sqrt{2\pi}}.
\end{eqnarray}
Let us denote 
\begin{equation}
F_4=2\sqrt{x}\;\hpg{2}{1}{1/2,\,1/4}{5/4}{\,x^2}, 
\qquad  C_4=\frac{\Gamma(1/4)^2}{\sqrt{2\pi}}.
\end{equation}
Formula (\ref{elin4a}) identifies $F_4$ as an elliptic integral $\int dx\big/\sqrt{x^3-x}$.
The action of the monodromy group on $F_4$ is generated by
$\sigma_0F_4=iF_4$, $\sigma_1F_4=C_4-F_4$.
The ``loop" paths to the other Riemann sheet of $\sqrt{x^3-x}$ around $x=0$, $x=1$ and $x=-1$
are represented by, respectively, $\sigma_0^2$, $\sigma_1$ and $\sigma_0^{-1}\sigma_1\sigma_0$.
In particular, the paths $\sigma_0^2\sigma_1$ and $\sigma_0\sigma_1\sigma_0$ are genuine loops on the Riemann surface of $\sqrt{x^3-x}$. Analytic continuation along them gives a pair of generating 
{\em periods} of the elliptic integral:
\begin{equation}
\sigma_0^2\sigma_1F_4=F_4+C_4, \qquad
\sigma_0\sigma_1\sigma_0F_4=F_4+iC_4.
\end{equation}
Analytic continuation along the Pochhammer path in (\ref{eq:pochha}) gives  
$F_4+(1-i)C_4$.

Similarly, for $(a,b,c)=(1/6,1/2,1/3)$ or $(1/3,1/3,1/3)$ we have the connection formulas
\begin{eqnarray*} 
3z^{1/6}\;\hpg{2}{1}{1/2,\,1/6}{7/6}{\,z}+(1-z)^{1/2}\;\hpg{2}{1}{1/2,\,5/6}{3/2}{1-z} &\equal&
\frac{3\;\Gamma(1/3)^3}{2^{7/3}\,\pi},\\
 z^{1/3}\,\hpg{2}{1}{1/3,\,2/3\,}{4/3}{\,z}+(1-z)^{1/3}\,\hpg{2}{1}{1/3,\,2/3}{4/3}{1-z}
&\equal& \frac{\Gamma(1/3)^3}{2\,\sqrt{3}\;\pi}.
\end{eqnarray*}
From here were can compute generating periods of $\int dx/\sqrt{x^3-1}$ and
$\int (X^3-1)^{-2/3}dX$. Recalling (\ref{elin6a}) and (\ref{elin3a}), we conclude that
the first integral has the periods
\[
\frac{\Gamma(1/3)^3}{2^{1/3}\,\pi}, \qquad (\omega+1)\frac{\Gamma(1/3)^3}{2^{1/3}\,\pi},
\qquad \mbox{given by the paths $\sigma_0^3\sigma_1$, $\sigma_0^{-1}\sigma_1\sigma_0$},
\]
while the second integral has the periods
\[
i\,\frac{\Gamma(1/3)^3}{\pi}, \qquad i\,(\omega+1)\frac{\Gamma(1/3)^3}{\pi},
\qquad \mbox{given by the paths $\sigma_0\sigma_1\sigma_0$, $\sigma_1^2\sigma_0$}.
\]

\bibliographystyle{alpha}
\bibliography{../hypergeometric}

\begin{thebibliography}{Gou81}

\bibitem[AK03]{andkitaev}
F.V. Andreev and A.V. Kitaev.
\newblock Some examples of ${RS}_3^2(3)$-transformations of ranks 5 and 6 as
  the higher order transformations for the hypergeometric function.
\newblock {\em Ramanujan J.}, 7(4):455--476, 2003.

\bibitem[Beu00]{beukrspriv}
F.~Beukers.
\newblock {P}rivate communication.
\newblock 2000.

\bibitem[Erd53]{bateman}
A.~Erd\'elyi, editor.
\newblock {\em Higher Transcendental Functions}, volume~I.
\newblock McGraw-Hill Book Company, New-York, 1953.

\bibitem[Gou81]{goursat}
E.~Goursat.
\newblock Sur l'\'equation diff\'erentielle lin\'eaire qui adment pour
  int\'egrale la s\'erie hyperg\'eom\'etrique.
\newblock {\em Ann. Sci. \'Ecole Noprm. Sup.(2)}, 10:S3--S142, 1881.

\bibitem[Sil86]{silverman1}
J.~H. Silverman.
\newblock {\em The arithmetic of elliptic curves}, volume 106 of {\em Grad.
  Texts in Math.}
\newblock Springer Verlaag, New-York, 1986.

\bibitem[Vid04]{algtgauss}
R.~Vid\=unas.
\newblock Algebraic transformations of {G}auss hypergeometric functions.
\newblock Available at {\sf http://arxiv.org/math.CA/0408269}, 2004.

\bibitem[Vid05]{thyperbolic}
R.~Vid\=unas.
\newblock Transformations of some {G}auss hypergeometric functions.
\newblock {\em Journ.~Comp. Applied Math.}, 178:473--487, 2005.
\newblock Available at {\sf http://arXiv.org/math.CA/0310436}.

\bibitem[Vid07]{degeneratehpg}
R.~Vid\=unas.
\newblock Degenerate {G}auss hypergeometric functions.
\newblock {\em Kyushu Journal of Mathematics}, 61:109--135, 2007.
\newblock Available at {\sf http://arxiv.org/math.CA/0407265}.

\bibitem[Vid08]{talggaus}
R.~Vid\=unas.
\newblock Transformations of algebraic {G}auss hypergeometric functions.
\newblock Available at {\sf http://arxiv.org/abs/0807.4808}, 2008.

\end{thebibliography}

\end{document}